\documentclass[11pt]{article}
\usepackage{pslatex}
\usepackage{fancyhdr}
\usepackage{graphicx}
\usepackage{geometry}
\usepackage{amssymb}
\usepackage{amsmath}
\usepackage{graphics}
\usepackage{epsfig}
\usepackage[dvips]{color}
\newtheorem{claim}{\bf \t}[part]

\newtheorem{Definition}{Definition}[part]

\newtheorem{Lemma}{Lemma}[part]
\newtheorem{Proposition}{Proposition}[part]
\newtheorem{Remark}{Remark}[part]
\newtheorem{Theorem}{Theorem}[part]

\numberwithin{Assumption}{section}
\numberwithin{Corollary}{section}
\numberwithin{Definition}{section}
\numberwithin{equation}{section} \numberwithin{Example}{section}
\numberwithin{Lemma}{section} \numberwithin{Proposition}{section}
\numberwithin{Remark}{section} \numberwithin{Theorem}{section}

\def\text#1{{\rm #1}}

\begin{document}
\title{\Large\bf Multivalued backward doubly stochastic differential equations with time delayed coefficients  \thanks{The work of Wen Lu is supported
partially by the National Natural Science Foundation of China
(61273128) and a Project of Shandong  Province Higher Educational
Science and Technology Program (J13LI06). The work of Yong Ren is
supported by the National Natural Science Foundation of
 China (10901003),
 the Distinguished Young Scholars of Anhui Province (1108085J08),
 the Key
 Project of Chinese Ministry of Education (211077) and the Anhui Provincial Natural
 Science Foundation (10040606Q30).  The work of Lanying Hu is
supported by the National Natural Science Foundation of
 China (11201004) }}

\author{\textbf{Wen Lu}$^1$ \footnote{e-mail: llcxw@163.com}\;\  \ \ \textbf{Yong
Ren}$^2$\footnote{Corresponding author. e-mail:
brightry@hotmail.com and renyong@126.com}\;\  \ \  \textbf{Lanying Hu}$^2$\footnote{e-mail: lanyinghu@126.com}
\\
 \small  1. School of Mathematics and Informational Science, Yantai University, Yantai 264005, China   \\
 \small 2. Department of Mathematics, Anhui Normal University, Wuhu
241000, China}\,
\date{}
\maketitle

\begin{abstract}
In this paper,  we deal with a class of multivalued backward doubly
 stochastic differential equations with time delayed coefficients.
 Based on a slight extension of the existence and uniqueness of solutions for
 backward doubly stochastic differential
equations with  time delayed coefficients, we establish the
existence and uniqueness of solutions for these equations by means
of Yosida approximation.
\end{abstract}

\vspace{.08in} \noindent \textbf{Keywords:}  backward doubly stochastic differential
equation;  time delayed coefficients; subdifferential operator;
Yosida approximation.

\vspace{.08in} \noindent \textbf{\bf MSC} 60H10, 60G40, 60H30

\section{Introduction }

Backward Stochastic Differential Equations (BSDEs in short) have
been first introduced in Pardoux and Peng \cite{PP90} in order to
give a probabilistic interpretation (Feynman-Kac formula) for the
solutions of semilinear parabolic PDEs.  In addition, in order to
give a probabilistic representation for a class of quasilinear
stochastic partial differential equations, Pardoux and Peng
\cite{Pardoux1994} introduced a new class of BSDEs, called as
backward doubly stochastic differential equations (BDSDEs in short).
This equation involves two different directions of stochastic
integrals and has also appeared as a powerful tool to give
probabilistic formulas for solutions of stochastic PDEs/PDIEs. One
can see Bally and Matoussi \cite{Bally2001}; Matoussi
\cite{Matoussi2002}; Zhang and Zhao \cite{Zhang2007}; Ren et al.
\cite{Ren2009} and the references therein.

On the other hand, BSDEs involving a subdifferential operator has
been treated by Pardoux and R\c{a}\u{s}canu \cite{PardouxRascanu98},
that they used to give a probabilistic representation for a class of
parabolic (and elliptic) variational inequalities. Furthermore,
Maticiuc and R\c{a}\u{s}canu \cite{MaticiucRascanu2010} gave a
probability interpretation of the viscosity solution of the
parabolic variational inequality (PVI in short) with a mixed
nonlinear multivalued Neumann--Dirichlet boundary condition.
Moreover, Boufoussi and Mrhardy \cite {Boufoussi2008} established
the existence result to stochastic viscosity solution for some
multivalued parabolic stochastic partial differential equation via
BDSDEs whose coefficient contains the subdifferential of a convex
function.

Recently, Delong and Imkeller \cite{DelongImkeller2010a} introduced
a class of BSDEs with time delayed generators of the form
\begin{eqnarray}\label{bdsde:00}
Y(t)=\xi+\int_t^T f(s, Y_s, Z_s)ds-\int_t^T Z(s)dW(s),
\end{eqnarray}
where the generator $f$ at time $s$ depends arbitrary on the past
values of a solution $(Y_s, Z_s)=(Y(s+u), Z(s+u)), -T\leq u\leq 0$.
They proved in \cite{DelongImkeller2010a} the  existence and
uniqueness of a solution for \eqref{bdsde:00}. Moreover, in Delong
and Imkeller \cite{Delong2010b}, they established the existence and
uniqueness as well as the Malliavin's differentiability of the
solution for BSDEs with time delayed generators driven by Brownian
motions and Poisson random measures.  Following this, Reis et al.
\cite{Reis2011} extended the results of  \cite{DelongImkeller2010a}
and \cite{Delong2010b} in $L^p$-spaces. For the applications of
BSDEs with time delayed coefficients in insurance and finance, one
can see Delong \cite{Delong2010}. Very recently, Diomande and
Maticiuc \cite{Diomande2013} established the existence and
uniqueness result result for multivalued BSDEs with time delayed
generators.

Besides, Lu and Ren \cite{LuRen2012} proved the existence and
uniqueness of the solutions for a class of BDSDEs with time delayed
coefficients under Lipschitz condition. Based on an extension of the
existence result of Lu and Ren \cite{LuRen2012}, the present paper
is to establish the existence and uniqueness of the solutions for
multivalued BDSDEs with time delayed coefficients under Lipschitz
condition.

The paper is organized as follows. In Section 2, we give some
preliminaries. Section 3 is concerned with BDSDEs with time delayed
coefficients. In Section 4,  we prove the existence and uniqueness
of the solution for multivalued BDSDEs with time delayed
coefficients.

\section{Notations, preliminaries and basic assumptions}

In this section, we provide some spaces and notations used in the
sequel. More precisely,  consider two mutually independent
$d$-dimensional Brownian motions $\{W_t, 0\leq t\leq T\}$ and
$\{B_t, 0\leq t\leq T\}$ defined on the probability spaces
$(\Omega_1, \mathcal {F}_1, \mathbb{P}_1)$ and $(\Omega_2, \mathcal
{F}_2, \mathbb{P}_2)$, respectively, where $T<\infty$ is a finite
time horizon. We denote
$$\mathcal {F}^{B}_{s,t}:=\sigma\{B_r-B_s, s\leq r\leq t\},\ \mathcal {F}^{W}_{t}:=\sigma\{W_r, 0\leq r\leq t\}. $$
Moreover, we define $\Omega=\Omega_1\times\Omega_2, \mathcal
{F}=\mathcal {F}_1\otimes\mathcal {F}_2$ and
$\mathbb{P}=\mathbb{P}_1\otimes \mathbb{P}_2$. We put
$$\mathcal {F}_t\triangleq \mathcal {F}^{W}_{t}\otimes\mathcal {F}^{B}_{s,t}\otimes\mathcal {N},$$
where $\mathcal{N}$ is the collection of $\mathbb{P}$-null sets of
$\mathcal {F}$. We use the usual Euclidian norm $|\cdot|$ in
$\mathbb{R}^k$ and $\mathbb{R}^{k\times d}$.

In what follows, we need the following spaces.
 \begin{itemize}
 \item $L^2_{-T}(\mathbb{R}^{k\times d})$: the space of measurable
functions $z: [-T,0]\rightarrow \mathbb{R}^{k\times d} $ such that
 $\int_{-T}^0|z(t)|^2dt<\infty.$

     \item  $L^{\infty}_{-T}(\mathbb{R}^k)$: the space of measurable
functions $y: [-T,0]\rightarrow \mathbb{R}^k$ such that
$\sup_{-T\leq t\leq 0}|y(t)|^2<\infty.$


\item $H^2_{T}(\mathbb{R}^{m})$: the space of
$\mathcal{F}$-predictable processes $\psi:
\Omega\times[0,T]\rightarrow \mathbb{R}^{m} $ such that
 $ E\int_{0}^T|\psi(t)|^2dt<\infty.$

\item $S^{2}_{T}(\mathbb{R}^k)$: the space of
$\mathcal{F}$-adapted, product measurable  processes $Y:
\Omega\times[0,T]\rightarrow \mathbb{R}^k $ such that
 $ E \left[\sup_{0\leq t\leq T}|Y(t)|^2 \right]<\infty.$
\end{itemize}

The spaces $H^2_{T}(\mathbb{R}^{k\times d})$ and
$S^{2}_{T}(\mathbb{R}^k)$ are endowed with the norm
$$ \|Z\|^2_{H_T^2}= E\int_{0}^T e^{\beta t}|Z(t)|^2dt \
\mbox{and} \  \|Y\|^2_{S_T^2}= E\left[\sup_{0\leq t\leq T} e^{\beta
t}|Y(t)|^2\right]$$ respectively with some $\beta>0$.
\bigskip

The purpose of the present paper is to consider the following
multivalued BDSDE with time delayed coefficients:
\begin{eqnarray}\label{mbdsde:1}
 \left\{ \begin{array}{l@{ }r}-dY(t)+\partial\varphi(Y(t))dt \ni f(t, Y(t), Z(t),
Y_t, Z_t)dt    +g(t, Y_t,
Z_t)dB(t)-Z(t)dW(t), \  \ 0\leq t\leq T,\\
Y_T=\xi,
\end{array}\right.
\end{eqnarray}
here the coefficients  $f$ and $g$ at time set can depend on the
past values of the solution denoted by $Y_s:=(Y(s+\theta))_{-T\leq
\theta\leq 0}$ and $Z_s:=(Z(s+\theta))_{-T\leq \theta\leq 0}$.
\begin{Remark}\label{remark:1} \rm
Throughout this paper, we always assume that $Y(t)=0$ and $Z(t)=0$
for $t<0$.
\end{Remark}

We mention that $\partial\varphi$ in Eq.\eqref{mbdsde:1} is the
subdifferential operator of the function $\varphi: \mathbb{R}^k
\rightarrow (-\infty, +\infty]$ which satisfies:
\\
(i) $\varphi$ is  proper ($\varphi\not\equiv \infty$), convex and
lower semicontinuous (l.s.c. for short),
\\
(ii) without loss generality, $\varphi(y)\geq\varphi(0)=0$, $\forall
y\in \mathbb{R}^k$.

\bigskip

For $\varphi$, let's define:
 \begin{itemize}
  \item ${\rm Dom}\varphi:=\{u\in \mathbb{R}^k: \varphi(u)<\infty\}$,

 \item $\partial\varphi(u):=\{u^*\in \mathbb{R}^k: \langle u^*,
v-u\rangle+\varphi(u)\leq\varphi(v), \forall v\in \mathbb{R}^k\}$,

 \item ${\rm Dom}(\partial\varphi):=\{u\in \mathbb{R}^k:
\partial\varphi(u)\neq\emptyset\}$,

 \item $(u, u^*)\in \partial\varphi \Longleftrightarrow u\in {\rm
Dom}(\partial\varphi), u^*\in \partial\varphi(u)$.
\end{itemize}
\begin{Remark} \rm
It is well known that the subdifferential operator $\partial\varphi$
is a maximal monotone operator, i.e., is maximal in the class of
operators which satisfy the condition
$$ \langle u^*-v^*,
u-v\rangle\geq 0,\; \forall (u, u^*), (v, v^*)\in \partial\varphi.$$
\end{Remark}
Now,  we give the following assumptions.
\begin{enumerate}
\item [(H1)] $E\left[|\xi|^2+\varphi(\xi)\right]<\infty$.

 \item [(H2)] The coefficients $f: \Omega\times [0,T]\times \mathbb{R}^{k}\times\mathbb{R}^{k\times d}\times
L^{\infty}_{-T}(\mathbb{R}^k)\times L^{2}_{-T}(\mathbb{R}^{k\times
d}) \rightarrow \mathbb{R}^k$ and $g: \Omega\times [0,T]\times
L^{\infty}_{-T}(\mathbb{R}^k)\times L^{2}_{-T}(\mathbb{R}^{k\times
d}) \rightarrow  \mathbb{R}^{k\times d}$ are product measurable,
$\mathcal{F}$-adapted and Lipschitz continuous in the sense that
there exist positive constant $K, L$ and $R$ such that, for a
non-random, finitely valued measure $\alpha$ supported on $[-T, 0)$
and for any $t\in[0,T]$, $(y^1,z^1), (y^2,z^2)\in \mathbb{R}^k\times
 \mathbb{R}^{k\times d}$, $(y^1_t,z^1_t), (y^2_t,z^2_t)\in
L^{\infty}_{-T}(\mathbb{R}^k)\times L^2_{-T}(\mathbb{R}^{k\times
d})$, $\mathbb{P}$-a.s.
\begin{enumerate}
\item [(i)] $\displaystyle|f(t,y^1,z^1,y^1_t,z^1_t)-f(t,y^2,z^2,y^1_t,z^1_t)|
 \leq K ( |y^1 -y^2 | + |z^1 -z^2 | )$;

\item [(ii)] $\displaystyle|f(t,y^1,z^1,y^1_t,z^1_t)-f(t,y^1,z^1,y^2_t,z^2_t)|^2  \leq
L\Big(\int_{-T}^0|y^1(t+\theta)-y^2(t+\theta)|^2\alpha(d\theta)$
$$+\int_{-T}^0 |z^1(t+\theta)-z^2(t+\theta) |^2\alpha(d\theta)\Big);$$

\item [(iii)] $\displaystyle|g(t,y^1,z^1,y^1_t,z^1_t)-g(t,y^2,z^2,y^1_t,z^1_t)|
 \leq R ( |y^1 -y^2 | + |z^1 -z^2 | )$;

\item [(iv)] $\displaystyle|g(t,y^1,z^1,y^1_t,z^1_t)-g(t,y^1,z^1,y^2_t,z^2_t)|^2  \leq
L\Big(\int_{-T}^0|y^1(t+\theta)-y^2(t+\theta)|^2\alpha(d\theta)$$$+\int_{-T}^0 |z^1(t+\theta)-z^2(t+\theta) |^2\alpha(d\theta)\Big).$$
\end{enumerate}
\item [(H3)] $\displaystyle
 \mathbb{E}\int_0^T|f(t,0,0,0,0)|^2dt<\infty,\ \  \mathbb{E}\int_0^T
|g(t,0,0) |^2dt<\infty.$
\item [(H4)] $\displaystyle f(t,\cdot,\cdot,\cdot,\cdot)=0,\ \ g(t,\cdot,\cdot)=0$ for $t<0$.
\end{enumerate}
\begin{Remark} \rm
We remark that, taking the measure $\alpha$ as Dirac measure
$\delta_{-r}$   with $r\in(0,T]$ or as Lebesgue measure,   the
coefficients could be of the form $k(t,y_t,z_t)=L\theta(t-r)$ or
$k(t,y_t,z_t)=L\int_0^t\theta(s)ds$ with $k=f,g$ and $\theta=y,z$.
\end{Remark}

We end this section by introduce the definition of the solution for
multivalued BDSDE \eqref{mbdsde:1}.
\begin{Definition}\label{def2:4} \rm
 The triple $(Y, Z, U)$ is a solution of multivalued
BDSDE \eqref{mbdsde:1} with subdifferential operator if
\begin{enumerate}
  \item [(i)]  $(Y,
Z, U)\in S^2_{T}(\mathbb{R}^{k})\times H^2_{T}(\mathbb{R}^{k\times
d})\times H^2_{T}(\mathbb{R}^{k})$,
  \item[(ii)] $\displaystyle E\int_0^Te^{\lambda t}\varphi(Y(t))dt<\infty$ ,
  \item [(iii)] $\displaystyle  (Y(t),U(t))\in \partial\varphi$,  $d\mathbb{P}\times dt$-a.e.
on $\Omega\times [0,T]$,
\item [(iv)] $\displaystyle  Y(t)+\int_t^TU(s)ds=\xi+\int_t^T
f(s,Y(s),Z(s),Y_s,Z_s)ds+\int_t^Tg(s,Y(s),Z(s),Y_s,Z_s)dB(s)-\int_t^TZ(s)dW(s)$,
\;$t\in[0,T]$.
\end{enumerate}
\end{Definition}
\section{BDSDEs with time delay coefficints}

In this part, we consider a class of BDSDEs with time delayed
coefficints as the form:
\begin{eqnarray}\label{bdsde:1}
Y(t)=\xi+\int_t^T f(s, Y(s), Z(s), Y_s, Z_s)ds +\int_t^T g(s, Y(s),
Z(s), Y_s, Z_s)dB(s)-\int_t^T Z(s)dW(s), 0\leq t\leq T.
\end{eqnarray}
We mention that the above equation is an extension of that of Lu and
Ren \cite{LuRen2012}, since the coefficients $f$ and $g$ can depend
on both the present and the past values of a solution $(Y, Z)$.

Now, we propose the definition of the solution for BDSDE
\eqref{bdsde:1}.
\begin{Definition} \label{definition:1} \rm
A solution to the BDSDE \eqref{bdsde:1} is  a pair of
$(Y,Z)=(Y(t),Z(t))_{0\leq t\leq T}$ satisfying that the BDSDE
\eqref{bdsde:1} such that $(Y,Z)\in S^{2}_{T}(\mathbb{R}^k)\times
H^2_{T}(\mathbb{R}^{k\times d})$.
\end{Definition}

Next, we list some results on BDSDE \eqref{bdsde:1}. Since their
proofs are similar to that of Lemma 3.1 and Theorem 3.2 of Lu and
Ren \cite{LuRen2012} with only a few slight changes, so we prefer to
omit them.

\begin{Lemma}\label{lemma:1}
Assume  $\mathbb{E}|\xi|^2<\infty$ and  the  assumptions
{\rm(H2)--(H4)} hold, and $(Y,Z)\in S^{2}_{T}(\mathbb{R}^k)\times
H^2_{T}(\mathbb{R}^{k\times d})$ be a solution of the
  BDSDED  {\rm\eqref{bdsde:1}}. If the Lipschitz constant $L$ of the coefficients $f$ and $g$ is small enough,
   then  there exist two positive constants $\beta$ and $C$ such that
\begin{eqnarray}
&&  \mathbb{E}\left[\sup_{0\leq s\leq T}{\rm e}^{\beta s}|Y(s)|^2+
\int_0^T {\rm e}^{\beta s}|Z(s)|^2ds\right] \nonumber\\&\leq&
C\mathbb{E}\left[{\rm e}^{\beta T}|\xi|^2 + \int_0^T {\rm e}^{\beta
s}|f(s,0,0,0,0)|^2ds+ \int_0^T {\rm e}^{\beta
s}|g(s,0,0,0,0)|^2ds\right].
 \end{eqnarray}
\end{Lemma}
Here and in the sequel, $C>0$ denotes a constant whose value may
change from one line to another.
\begin{Theorem}\label{theorem:1}
Assume $\mathbb{E}|\xi|^2<\infty$ and  the assumptions
{\rm(H2)--(H4)} hold. If the Lipschitz constant $L$ of the
coefficients $f$ and $g$ is small enough , then the BDSDE  {\rm
\eqref{bdsde:1}} has a unique solution.
\end{Theorem}

\section{Existence and uniqueness of the solution}

This section is devoted to the study of the existence and uniqueness
result of  multivalued BDSDE  \eqref{mbdsde:1}.

The main result of this section is the following
\begin{Theorem}\label{theorem4:1}
Assume that the assumptions {\rm(H1)--(H4)} hold. If the Lipschitz
constant $L$ of the  coefficients $f$ and $g$ is small enough , then
there exists a unique solution of BDSDE  {\rm \eqref{mbdsde:1}}.
\end{Theorem}

We mention that our proof is based on the Yosida approximations.
 First of all, let's introduce
an approximation of the function $\varphi$ by a convex
$C^1$--function $\varphi_\epsilon$, $\epsilon>0$, defined by
\begin{eqnarray}\label{approxi:2}
\varphi_\epsilon(u)&=&\inf\left\{\frac{1}{2\epsilon}|u-v|^2+\varphi(v):
v\in \mathbb{R}^k\right\}\nonumber\\&=&\frac{1}{2\epsilon}|u-J_\epsilon
(u)|^2+\varphi(J_\epsilon u),
\end{eqnarray}
where $J_\epsilon (u)=(I+ \epsilon\partial\varphi)^{-1}(u)$. For
convenience, we illustrate some properties of this approximation,
for more details, one can see Br\'{e}zis \cite{Brezis1973}.
\begin{Proposition}\label{prop4:1}
 For all $\epsilon, \delta>0$, $u, v\in
\mathbb{R}^k$, it holds that
\begin{enumerate}
  \item [\rm(i)] $\varphi_\epsilon$ is a convex function with the
gradient  being a
Lipschitz function;
  \item [\rm(ii)] $\displaystyle \varphi_\epsilon(u)\leq \varphi(u)$;
  \item [\rm(iii)] $\displaystyle \nabla
\varphi_\epsilon(u)=\partial\varphi_\epsilon(u)=\frac{u-J_\epsilon
(u)}{\epsilon}\in\partial\varphi(J_\epsilon (u))$;
  \item [\rm(iv)] $\displaystyle |J_\epsilon (u)-J_\epsilon (v)|\leq  |u-v|$;
  \item [\rm(v)] $\displaystyle 0\leq \varphi_\epsilon(u)\leq
\langle\nabla\varphi_\epsilon(u), u\rangle$;
  \item [\rm(vi)] $\displaystyle\langle\nabla\varphi_\epsilon(u)-\nabla\varphi_\delta(v),
u-v\rangle\geq-(\epsilon+\delta)\langle\nabla\varphi_\epsilon(u),
\nabla\varphi_\delta(v)\rangle$.
\end{enumerate}
\end{Proposition}
Let us consider the approximating equation as the form:
\begin{eqnarray}\label{bdsde:2}
Y^{\epsilon}_t+\int_t^T\nabla\varphi_{\epsilon}(Y^{\epsilon}(s))ds&=&
\xi+\int_t^T
f(s,Y^{\epsilon}(s),Z^{\epsilon}(s),Y^{\epsilon}_s,Z^{\epsilon}_s)ds\nonumber\\&&+\int_t^T
g(s,Y^{\epsilon}(s),Z^{\epsilon}(s),Y^{\epsilon}_s,Z^{\epsilon}_s)dB(s)-\int_t^TZ^{\epsilon}(s)dW(s),
 t\in[0,T],
\end{eqnarray}
where $\xi, f$ and $g$ satisfy assumptions (H1)--(H4). Since
$\nabla\varphi_{\epsilon}$ is Lipschitz continuous, we know that
from
 Theorem \ref{theorem:1}, for  sufficiently small $L$ and $R$, BDSDE \eqref{bdsde:2} has a
 unique solution $(Y^{\epsilon}, Z^{\epsilon})\in S^{2}_{T}(\mathbb{R}^k)\times
H^2_{T}(\mathbb{R}^{k\times d})$.
\begin{Lemma}\label{lemma4:1}
Assume the assumptions {\rm(H1)--(H4)} hold. If the Lipschitz
coefficients $L$ and $R$ are small enough, then it holds that
\begin{eqnarray}
\mathbb{E}\left[\sup_{t\in [0,T]}e^{\beta
t}|Y^\epsilon(t)|^2+\int_0^Te^{\beta s} |Z^\epsilon(s)
|^2ds\right]\leq CM_1,
\end{eqnarray}
 where $M_1:=\mathbb{E}\left[ e^{\beta T}|\xi|^2+\int_0^Te^{\beta s} |f(s,0,0,0,0)
|^2ds+\int_0^Te^{\beta s} |g(s,0,0,0,0) |^2ds\right].$
\end{Lemma}
\noindent {\bf Proof.} For any $\beta>0$, applying It\^{o}'s formula to $e^{\beta
t}|Y^\epsilon(t)|^2$ yields that

$$e^{\beta t}|Y^\epsilon(t)|^2+\beta\int_t^Te^{\beta
s}|Y^\epsilon(s)|^2ds+2\int_t^Te^{\beta s}\langle
Y^\epsilon(s), \nabla\varphi_{\epsilon}(Y^\epsilon(s))\rangle ds+
\int_t^Te^{\beta s} |Z^\epsilon(s) |^2ds$$
\begin{eqnarray}\label{lm1:1}
&=&e^{\beta T}|\xi|^2+2\int_t^Te^{\beta s}\langle Y^\epsilon(s),
f(s,Y^\epsilon(s),Z^\epsilon(s),Y^\epsilon_s, Z^\epsilon_s)\rangle
ds+\int_t^Te^{\beta s}
|g(s,Y^\epsilon(s),Z^\epsilon(s),Y^\epsilon_s,
Z^\epsilon_s)|^2ds\nonumber\\&&+2\int_t^Te^{\beta s}\langle
Y^\epsilon(s), g(s,Y^\epsilon(s),Z^\epsilon(s), Y^\epsilon_s,
Z^\epsilon_s) dB(s)\rangle -2\int_t^Te^{\beta s}\langle
Y^\epsilon(s), Z^\epsilon(s)dW(s)\rangle.
\end{eqnarray}
By Young's inequality and (H2), we have
\begin{eqnarray}\label{lm1:2}
&& 2\int_t^Te^{\beta s}\langle Y^\epsilon(s),
f(s,Y^\epsilon(s),Z^\epsilon(s),Y^\epsilon_s, Z^\epsilon_s)\rangle
ds\nonumber
\\&\leq& \gamma\int_t^Te^{\beta
s}|Y^\epsilon(s)|^2ds +\frac{1}{\gamma}\int_t^Te^{\beta
s}|f(s,Y^\epsilon(s),Z^\epsilon(s),Y^\epsilon_s,
Z^\epsilon_s)|^2ds\nonumber
\\&\leq& \gamma\int_t^Te^{\beta
s}|Y^\epsilon(s)|^2ds +\frac{1}{\gamma}\int_t^Te^{\beta
s}|f(s,Y^\epsilon(s),Z^\epsilon(s),Y^\epsilon_s,
Z^\epsilon_s)|^2ds\nonumber
\\&\leq& \gamma\int_t^Te^{\beta
s}|Y^\epsilon(s)|^2ds +\frac{3}{\gamma}\int_t^Te^{\beta
s}|f(s,0,0,0, 0)|^2ds
+\frac{6K^2}{\gamma}\int_t^Te^{\beta s}(|Y^\epsilon(s)|^2+
|Z^\epsilon(s) |^2)ds\nonumber
\\&&
+\frac{3L}{\gamma}\int_t^Te^{\beta
s}\int_{-T}^0(|Y^\epsilon(s+\theta)|^2+ |Z^\epsilon(s+\theta)
|^2)\alpha(d\theta)ds
 \end{eqnarray}
and
\begin{eqnarray}\label{lm1:3}
&&\int_t^Te^{\beta s} |g(s,Y^\epsilon(s),Z^\epsilon(s),Y^\epsilon_s,
Z^\epsilon_s)|^2ds \leq 6R^2\int_t^Te^{\beta
s}[|Y^\epsilon(s)|^2+|Z^\epsilon(s)|^2]ds\nonumber\\&&+3\int_t^Te^{\beta
s} |g(s,0,0,0,0)|^2ds +3L\int_t^Te^{\beta
s}\int_{-T}^0(|Y^\epsilon(s+\theta)|^2+ |Z^\epsilon(s+\theta)
|^2)\alpha(d\theta)ds.
 \end{eqnarray}
By a change of integration order argument and Remark \ref{remark:1},
we obtain
\begin{eqnarray}\label{lm1:4}
&&\int_0^T e^{\beta
s}\int_{-T}^0|Y(s+\theta)|^2\alpha(d\theta)ds=\int_{-T}^0\int_0^T
e^{\beta s}|Y(s+\theta)|^2ds\alpha(d\theta)\nonumber \\&=&
\int_{-T}^0{\rm e}^{-\beta \theta}\int_\theta^{T+\theta} e^{\beta
t}|Y(t)|^2dt\alpha(d\theta) \leq
\min\left\{T\widetilde{\alpha}\sup_{0\leq t\leq T} e^{\beta
t}|Y(t)|^2, \; \widetilde{\alpha}\int_0^{T} e^{\beta
t}|Y(t)|^2dt\right\}.
\end{eqnarray}
and
\begin{eqnarray}\label{lm1:5}
\int_0^T e^{\beta s}\int_{-T}^0 |Z(s+\theta)
|^2\alpha(d\theta)ds=\int_{-T}^0{\rm e}^{-\beta
\theta}\int_\theta^{T+\theta} e^{\beta t}|Z(t)|^2dt\alpha(d\theta)
\leq  \widetilde{\alpha} \int_0^{T} e^{\beta t} |Z(t) |^2dt,
\end{eqnarray}
 where $\widetilde{\alpha}=\int_{-T}^0 e^{-\beta \theta}\alpha(d\theta)$.

Combining \eqref{lm1:2}-\eqref{lm1:5} together with (v) of
Proposition \ref{prop4:1}, \eqref{lm1:1} becomes
\begin{eqnarray}\label{lm1:6}
&&e^{\beta t}|Y^\epsilon(t)|^2+\beta\int_t^Te^{\beta
s}|Y^\epsilon(s)|^2ds+\int_t^Te^{\beta s} |Z^\epsilon(s) |^2ds
\nonumber
\\&\leq& e^{\beta T}|\xi|^2+\gamma\int_t^Te^{\beta
s}|Y^\epsilon(s)|^2ds+\frac{6K^2}{\gamma}\int_t^Te^{\beta
s}(|Y^\epsilon(s)|^2+ |Z^\epsilon(s) |^2)ds\nonumber
\\&&
+\frac{3}{\gamma}\int_t^Te^{\beta s}|f(s,0,0,0,
0)|^2ds+3\int_t^Te^{\beta s} |g(s,0,0,0, 0) |^2ds \nonumber
\\&&
+ 6R^2  \int_t^{T}e^{\beta s}
(|Y^\epsilon(s)|^2+|Z^\epsilon(s)|^2)ds
+\left(\frac{3L\widetilde{\alpha}}{\gamma}+3L\widetilde{\alpha}\right)
\int_0^{T}e^{\beta s}(|Y^\epsilon(s)|^2+|Z^\epsilon(s)|^2)ds
\nonumber
\\&&
+2\int_t^Te^{\beta s}\langle Y^\epsilon(s), g(s, Y^\epsilon_s,
Z^\epsilon_s) dB(s)\rangle -2\int_t^Te^{\beta s}\langle
Y^\epsilon(s), Z^\epsilon(s)dW(s)\rangle.
\end{eqnarray}
For $t=0$, taking expectation on both sides of above gives
\begin{eqnarray}\label{lm1:7}
&& \mathbb{E}|Y^\epsilon(0)|^2+K_1\mathbb{E}\int_0^Te^{\beta
s}|Y^\epsilon(s)|^2ds +K_2\mathbb{E}\int_0^Te^{\beta s}
|Z^\epsilon(s) |^2ds \nonumber
\\&\leq& \mathbb{E}[e^{\beta T}|\xi|^2]+\frac{3}{\gamma}\mathbb{E}\int_0^Te^{\beta s}|f(s,0,0,0,
0)|^2ds+3\mathbb{E}\int_0^Te^{\beta s} |g(s,0,0,0,0)|^2ds.
\end{eqnarray}
where
$K_1:=\beta-\gamma-\frac{6K^2}{\gamma}-\frac{3L\widetilde{\alpha}}{\gamma}-
3L\widetilde{\alpha}-6R^2$,
$K_2:=1-\frac{6K^2}{\gamma}-\frac{3L\widetilde{\alpha}}{\gamma}
-3L\widetilde{\alpha}-6R^2$.

For sufficiently small $L$ and $R$, choosing $\beta>0$, $\gamma>0$
such that $K_1>0$ and $K_2>0$, by \eqref{lm1:7},  we obtain that
there exists a constant $C>0$ depending on $\beta, \gamma, K, L, R$
and $\widetilde{\alpha}$ such that
\begin{eqnarray}\label{lm1:8}
&&\mathbb{E}\int_0^Te^{\beta
s}|Y^\epsilon(s)|^2ds+\mathbb{E}\int_0^Te^{\beta s} |Z^\epsilon(s)
|^2ds \nonumber\\&\leq&  C\left\{ \mathbb{E}[e^{\beta
T}|\xi|^2]+\mathbb{E}\int_0^Te^{\beta s}|f(s,0,0,0, 0)|^2ds
+\mathbb{E}\int_0^Te^{\beta s} |g(s,0,0,0,0)|^2ds\right\}.
\end{eqnarray}
On the other hand, for  $\beta$ and $\gamma$ choosing above, from
\eqref{lm1:6}, we have
\begin{eqnarray}\label{lm1:9}
\sup_{0\leq t\leq T}e^{\beta t}|Y^\epsilon(t)|^2 &\leq& e^{\beta
T}|\xi|^2+\frac{3}{\gamma}\int_t^Te^{\beta s}|f(s,0,0,0,
0)|^2ds+3\int_t^Te^{\beta s} |g(s,0,0,0,0)|^2ds
 \nonumber
\\&&
+2\sup_{0\leq t\leq T}\left|\int_t^Te^{\beta s}\langle Y^\epsilon(s),
g(s,Y^\epsilon(s), Z^\epsilon(s), Y^\epsilon_s, Z^\epsilon_s)
dB(s)\rangle\right| \nonumber
\\&&
+2\sup_{0\leq t\leq
T}\left|\int_t^Te^{\beta s}\langle Y^\epsilon(s),
Z^\epsilon(s)dW(s)\rangle\right|.
\end{eqnarray}
By the Burkholder--Davis--Gundy inequality and Young's
 inequality, together with  \eqref{lm1:2}-\eqref{lm1:3} and (H2),
 there exists a constant $d_1>0$ such that
 \begin{eqnarray}\label{lm1:10}
&&  2 \mathbb{E}\left[\sup_{0\leq t\leq T}\left|\int_t^T  e^{\beta
s}\langle Y^\epsilon(s), g(s, Y^\epsilon(s), Z^\epsilon(s),
Y^\epsilon_s, Z^\epsilon_s)dB(s)\rangle \right|\right]\nonumber
\\&\leq&  d_1\Bigg[ \rho_1  \mathbb{E}\left(\sup_{0\leq t\leq T} e^{\beta
t}|Y^\epsilon(t)|^2\right) +\frac{6R^2+3L\widetilde{\alpha}}{\rho_1}
\mathbb{E} \int_0^T e^{\beta s} |Z^\epsilon(s) |^2ds\nonumber
\\&&
\quad+\frac{3}{\rho_1} \mathbb{E}\int_0^T e^{\beta s}|g(s,0,0,0,0)
|^2ds \Bigg].
 \end{eqnarray}
 Similarly, there exists  a constant $d_2>0$ such that
 \begin{eqnarray}\label{lm1:11}
&&  2 \mathbb{E}\left[\sup_{0\leq t\leq T}\left|\int_t^T  e^{\beta
s}\langle Y^\epsilon(s), Z^\epsilon(s)dW(s)\rangle
\right|\right]\nonumber
\\&\leq&  d_2\left[ \rho_2  \mathbb{E} \left(\sup_{0\leq t\leq T} e^{\beta
t}|Y^\epsilon(t)|^2\right) +\frac{1}{\rho_2} \mathbb{E}\int_0^T e^{\beta
s} |Z^\epsilon(s) |^2ds\right],
 \end{eqnarray}
 where $\rho_1$ and $\rho_2$ are two positive constants.

Then, choosing $\rho_1=\frac{1}{3d_1}$ and $\rho_2=\frac{1}{3d_2}$,
for sufficiently small $L>0$ and $R>0$,
there exists a constant $C>0$ depending on $\beta, \gamma, L, K$,
$\widetilde{\alpha}$, $d_1$ and $d_2$ such that
\begin{eqnarray*}
&&\mathbb{E}\left[\sup_{t\in [0,T]}e^{\beta
t}|Y^\epsilon(t)|^2+\int_0^Te^{\beta s} |Z^\epsilon(s) |^2ds\right]
\\&\leq& C\mathbb{E}\left[ e^{\beta T}|\xi|^2+\int_0^Te^{\beta
s} |f(s,0,0,0,0) |^2ds+\int_0^Te^{\beta s} |g(s,0,0,0,0)|^2ds\right].
\end{eqnarray*}
The lemma is proved.

\begin{Lemma}\label{lemma4:2}
 Assume that the conditions of Lemma \ref{lemma4:1} hold. Then, for all $0\leq t\leq T$, it holds that
\begin{enumerate}
  \item [\rm(i)] $\displaystyle \mathbb{E}\int_0^Te^{\beta s}|\nabla
\varphi_\epsilon(Y^\epsilon(s))|^2ds\leq CM_2$,
  \item [\rm(ii)] $\displaystyle \mathbb{E}\left[e^{\beta t}  \varphi
(J_\epsilon(Y^\epsilon(t)))+E\int_0^Te^{\beta s}\varphi
(J_\epsilon(Y^\epsilon(s)))ds\right]\leq CM_2$,
  \item [\rm(iii)]
 $\displaystyle \mathbb{E}[e^{\beta t}|Y^\epsilon(t)-
 J_\epsilon(Y^\epsilon(t))|^2ds]\leq C\epsilon M_2$,
\end{enumerate}
 where $M_2:=M_1+\mathbb{E}[e^{\beta
T}\varphi(\xi)]$.
\end{Lemma}
\noindent {\bf Proof.}  The stochastic subdifferential inequality in Pardoux
and R\u{a}\c{s}canu \cite {PardouxRascanu98} gives that
\begin{eqnarray*}
e^{\beta T}\varphi_{\epsilon}(\xi)  \geq  e^{\beta
t}\varphi_{\epsilon}(Y^{\epsilon}(t))+\int_t^Te^{\beta
s}\langle\nabla\varphi_{\epsilon}(Y^{\epsilon}(s)),
dY^{\epsilon}(s)\rangle
+\int_t^T\varphi_{\epsilon}(Y^{\epsilon}(s))d(e^{\beta s}).
\end{eqnarray*}
Therefore,
\begin{eqnarray}\label{lm1:15}
&& e^{\beta
t}\varphi_{\epsilon}(Y^{\epsilon}(t))+\beta\int_t^Te^{\beta
s}\varphi_{\epsilon}(Y^{\epsilon}(s))ds+\int_t^Te^{\beta
s}|\nabla\varphi_{\epsilon}(Y^{\epsilon}(s))|^2ds \nonumber\\&\leq&
e^{\beta T}\varphi_{\epsilon}(\xi) + \int_t^Te^{\beta s}\langle
\nabla\varphi_{\epsilon}(Y^{\epsilon}(s)),
f(s,Y^\epsilon(s),Z^\epsilon(s),Y^\epsilon_s, Z^\epsilon_s)\rangle
ds\nonumber\\&& + \int_t^Te^{\beta s}\langle
\nabla\varphi_{\epsilon}(Y^{\epsilon}(s)), g(s,Y^\epsilon(s),
Z^\epsilon(s), Y^\epsilon_s, Z^\epsilon_s) dB(s)\rangle
\nonumber\\&&-\int_t^Te^{\beta s}\langle
\nabla\varphi_{\epsilon}(Y^{\epsilon}(s)),
Z^\epsilon(s)dW(s)\rangle.
\end{eqnarray}
Since
\begin{eqnarray*}
&&  \int_t^Te^{\beta s}\langle
\nabla\varphi_{\epsilon}(Y^{\epsilon}(s)),
f(s,Y^\epsilon(s),Z^\epsilon(s),Y^\epsilon_s, Z^\epsilon_s)\rangle
ds\nonumber\\&\leq& \frac{1}{2}\int_t^Te^{\beta s}
|\nabla\varphi_{\epsilon}(Y^{\epsilon}(s))|^2ds+\frac{3}{2}\int_t^Te^{\beta
s} |f(s, 0, 0, 0, 0)|^2ds
\\&&
+3K^2\int_t^T(| Y^{\epsilon}(s)|^2+ | Y^{\epsilon}(s) |^2)ds
+\frac{3L\widetilde{\alpha}}{2} \int_0^T e^{\beta s}
(|Y^{\epsilon}(s)|^2+|Z^{\epsilon}(s)|^2)ds.
\end{eqnarray*}
Then, by Lemma \ref{lemma4:1} and the nonnegative property of
$\varphi_{\epsilon}(y)$, (i) is hold for sufficiently small $L$ and
$R$.

From \eqref{lm1:15}, for sufficiently small $L$ and $R$, we get
\begin{eqnarray*}
\mathbb{E}[e^{\beta t}\varphi_{\epsilon}(Y^{\epsilon}(t))]+
\mathbb{E}\int_t^Te^{\beta
s}\varphi_{\epsilon}(Y^{\epsilon}(s))ds\leq CM_2.
\end{eqnarray*}
since $\varphi(J_\epsilon (y))\leq\varphi_{\epsilon}(y)$ (see
\eqref{approxi:2}), it follows that
\begin{eqnarray*}
\mathbb{E}\left[e^{\beta
t}\varphi(J_{\epsilon}(Y^{\epsilon}(t)))\right]+\mathbb{E}\int_t^Te^{\beta
s}\varphi(J_{\epsilon}(Y^{\epsilon}(s)))ds\leq C M_2.
\end{eqnarray*}
Moreover, since
$$\frac{1}{2\epsilon}e^{\beta t}|Y^{\epsilon}(t)-J_{\epsilon}(Y^{\epsilon}(t))|^2
\leq e^{\beta t}\varphi_{\epsilon}(Y^{\epsilon}(t)),$$ we then have
$$\mathbb{E}\left[e^{\beta t}|Y^{\epsilon}(t)-J_{\epsilon}(Y^{\epsilon}(t))|^2\right]
\leq C \epsilon M_2.$$ The proof is complete.

\begin{Lemma}\label{lemma4:3}
 Assume that the conditions of Lemma \ref{lemma4:1} hold. Then, it holds that
  $$\mathbb{E}\left[\sup_{t\in [0,T]}e^{\beta t}|Y^\epsilon(t)-Y^\delta(t)|^2+\int_0^Te^{\beta
s} |Z^\epsilon(s)-Z^\delta(s) |^2ds\right]\leq
C(\epsilon+\delta)M_2.$$
\end{Lemma}
\noindent {\bf Proof.}  Applying It\^{o}'s formula to $e^{\beta
t}|Y^\epsilon(t)-Y^\delta(t)|^2$ yields
\begin{eqnarray}\label{lm1:18}
&&e^{\beta t}|Y^\epsilon(t)-Y^\delta(t)|^2+\beta\int_t^Te^{\beta
s}|Y^\epsilon(s)-Y^\delta(s)|^2ds+\int_t^Te^{\beta s}
|Z^\epsilon(s)-Z^\delta(s) |^2ds\nonumber\\&&+2\int_t^Te^{\beta
s}\langle Y^\epsilon(s)-Y^\delta(s),
 \nabla\varphi_{\epsilon}(Y^\epsilon(s))-\nabla\varphi_{\delta}(Y^\delta(s)) \rangle
ds\nonumber
\\&=& 2\int_t^Te^{\beta s}\langle
Y^\epsilon(s)-Y^\delta(s),
 f(s,Y^\epsilon(s),Z^\epsilon(s),Y^\epsilon_s,
Z^\epsilon_s)-f(s,Y^\delta(s),Z^\delta(s),Y^\delta_s, Z^\delta_s)
\rangle ds\nonumber\\&&+\int_t^Te^{\beta s}
|g(s,Y^\epsilon(s),Z^\epsilon(s),Y^\epsilon_s,
Z^\epsilon_s)-g(s,Y^\delta(s),Z^\delta(s),Y^\delta_s,
Z^\delta_s)|^2ds \nonumber
\\&&
+2\int_t^Te^{\beta s}\langle Y^\epsilon(s)-Y^\delta(s),
(g(s,Y^\epsilon(s),Z^\epsilon(s), Y^\epsilon_s,
Z^\epsilon_s)-g(s,Y^\delta(s),Z^\delta(s),Y^\delta_s,
Z^\delta_s))dB(s)\rangle\nonumber
\\&&
-2\int_t^Te^{\beta s}\langle Y^\epsilon(s)-Y^\delta(s),
(Z^\epsilon(s)-Z^\delta(s))dW(s)\rangle.
\end{eqnarray}
Since
\begin{eqnarray*}
\langle Y^\epsilon(s)-Y^\delta(s),
\nabla\varphi_{\epsilon}(Y^\epsilon(s))-\nabla\varphi_{\delta}(Y^\delta(s))\rangle\geq
-(\epsilon+\delta)|\nabla\varphi_{\epsilon}(Y^\epsilon(s))||\nabla\varphi_{\delta}(Y^\delta(s))|,
\end{eqnarray*}
by Young's inequality and the assumptions of $f$ and $g$, we get
\begin{eqnarray}\label{lm1:20}
&& 2\int_t^Te^{\beta s}\langle Y^\epsilon(s)-Y^\delta(s),
f(s,Y^\epsilon(s),Z^\epsilon(s),Y^\epsilon_s,
Z^\epsilon_s)-f(s,Y^\delta(s),Z^\delta(s),Y^\delta_s,
Z^\delta_s)\rangle ds \nonumber\\&\leq& \gamma\int_t^Te^{\beta
s}|Y^\epsilon(s)-Y^\delta(s)|^2ds
+\frac{4K^2}{\gamma}\int_t^Te^{\beta
s}(|Y^\epsilon(s)-Y^\delta(s)|^2+ |Z^\epsilon(s)-Z^\delta(s) |^2)ds
\nonumber\\&&+\frac{2L}{\gamma}\int_t^Te^{\beta
s}\left[\int_{-T}^0(|Y^\epsilon(s+\theta)-Y^\delta(s+\theta)|^2 +
|Z^\epsilon(s+\theta)-Z^\delta(s+\theta) |^2)\alpha(d\theta)\right]ds
\nonumber\\&\leq&  \gamma\int_t^Te^{\beta
s}|Y^\epsilon(s)-Y^\delta(s)|^2ds
+\frac{4K^2}{\gamma}\int_t^Te^{\beta
s}(|Y^\epsilon(s)-Y^\delta(s)|^2+ |Z^\epsilon(s)-Z^\delta(s) |^2)ds
\nonumber\\&&+\frac{2L\widetilde{\alpha}}{\gamma}\int_0^Te^{\beta
s}(|Y^\epsilon(s)-Y^\delta(s)|^2+ |Z^\epsilon(s)-Z^\delta(s) |^2)ds
\end{eqnarray}
and
\begin{eqnarray}\label{lm1:21}
&& \int_t^Te^{\beta s} |g(s,Y^\epsilon(s),Z^\epsilon(s),
Y^\epsilon_s, Z^\epsilon_s)-g(s,Y^\delta(s),Z^\delta(s),Y^\delta_s,
Z^\delta_s)|^2ds\nonumber\\&\leq&
 4R^2\int_t^Te^{\beta
s}(|Y^\epsilon(s)-Y^\delta(s)|^2+ |Z^\epsilon(s)-Z^\delta(s) |^2)ds
\nonumber\\&&+2L\widetilde{\alpha}\int_0^Te^{\beta
s}(|Y^\epsilon(s)-Y^\delta(s)|^2+ |Z^\epsilon(s)-Z^\delta(s) |^2)ds.
\end{eqnarray}
Combining \eqref{lm1:18}-\eqref{lm1:21}, we obtain
\begin{eqnarray}\label{lm1:22}
&& |Y^\epsilon(0)-Y^\delta(0)|^2+K_3\int_0^Te^{\beta
s}|Y^\epsilon(s)-Y^\delta(s)|^2ds +K_4\int_0^Te^{\beta s}
|Z^\epsilon(s)-Z^\delta(s) |^2ds \nonumber\\&\leq&
2(\epsilon+\delta)\int_t^Te^{\beta
s}|\nabla\varphi_{\epsilon}(Y^\epsilon(s))||\nabla\varphi_{\delta}(Y^\delta(s))|ds.
\end{eqnarray}
where $K_3:=\beta-\gamma-\frac{4K^2}{\gamma}
-\frac{2L\widetilde{\alpha}}{\gamma}-2L\widetilde{\alpha}-4R^2$,
$K_4:=1-\frac{4K^2}{\gamma}
-\frac{2L\widetilde{\alpha}}{\gamma}-2L\widetilde{\alpha}-4R^2$. For
 sufficiently small $L$ and $R$, choosing $\beta$, $\gamma>0$ such that
$K_3>0$ and $K_4>0$, by (i) of Lemma \ref{lemma4:2}, we have
\begin{eqnarray}\label{lm1:23}
\mathbb{E}\int_0^Te^{\beta s}|Y^\epsilon(s)-Y^\delta(s)|^2ds
+\mathbb{E}\int_0^Te^{\beta s} |Z^\epsilon(s)-Z^\delta(s) |^2ds \leq
C(\epsilon+\delta)M_2.
\end{eqnarray}
Therefore, as the same procedure as \eqref{lm1:10}--\eqref{lm1:11},
we can get the desired result from the Burkholder--Davis--Gundy
inequality and \eqref{lm1:23}.

Now, let's give a proof of Theorem \ref{theorem4:1}.

\noindent {\bf Proof.}  {\bf Existence.}  Lemma \ref{lemma4:3}
implies that there exist $Y\in S^{2}_{T}(\mathbb{R}^k)$ and $Z\in
H^2_{T}(\mathbb{R}^{k\times d})$ such that
$$\lim_{\epsilon\rightarrow0}(Y^\epsilon, Z^\epsilon)=(Y, Z).$$ Then
Lemma \ref{lemma4:2} shows that
\begin{eqnarray}
\lim_{\epsilon\rightarrow0}J_{\epsilon}(Y^\epsilon)=y  \ \ {\rm in}\
\  H^2_{T}(\mathbb{R}^{k})\nonumber
\end{eqnarray}
and
\begin{eqnarray}
\lim_{\epsilon\rightarrow0}E[e^{\beta
t}|J_{\epsilon}(Y^\epsilon(t))-y(t)|^2]=0,\ \ 0\leq t\leq
T.\nonumber
\end{eqnarray}
Moreover, Fatou's lemma, Lemma \ref{lemma4:2}, Proposition
\ref{prop4:1} and the lower semicontinuity of $\varphi$ shows that
(ii) of Definition \ref{def2:4} is satisfied.

In addition, (i) of Lemma \ref{lemma4:2} shows that
$U^\epsilon(t):=\nabla\varphi_{\epsilon}(Y^\epsilon(t))$ are bounded
in the space $H^2_{T}(\mathbb{R}^{k})$, so there exists a
subsequence $\epsilon_n\rightarrow 0$ such that
\begin{eqnarray}
U^{\epsilon_n} \rightarrow U, \ \ {\rm weakly \ in}\ \
H^2_{T}(\mathbb{R}^{k}).\nonumber
\end{eqnarray}
Furthermore, we have
\begin{eqnarray}
\mathbb{E}\int_0^T|U(s)|^2ds\leq\liminf_{n\rightarrow\infty}\mathbb{E}\int_0^T|U^{\epsilon_n}(s)|^2ds\leq
CM_2.\nonumber
\end{eqnarray}
In virtue of (H2), by passing limit in BDSDE \eqref{bdsde:2}, we
deduce that the triple $(Y,Z,U)$ satisfies (iv) of Definition
\ref{def2:4}.

Finally, let us show (iii) of  Definition \ref{def2:4} is satisfied.
 Since
 $U^\epsilon(t)\in\partial\varphi(J_{\epsilon}(Y^{\epsilon}(t)))$,
 $t\in[0,T]$, it follows that, for all $V\in
 H^2_{T}(\mathbb{R}^{k})$,
\begin{eqnarray}
e^{\beta t}\langle U^\epsilon(t),
V(t)-J_{\epsilon}(Y^{\epsilon}(t))\rangle+e^{\beta
t}\varphi(J_{\epsilon}(Y^{\epsilon}(t)))\leq e^{\beta
t}\varphi(V(t)), d\mathbb{P}\times dt-a.e.\nonumber
\end{eqnarray}
Taking the $\liminf$ in the probability in the above inequality,
(iii) of  Definition \ref{def2:4} holds.

{\bf Uniqueness.}    Let $(Y^i(t), Z^i(t), U^i(t))$, $i=1, 2$ be two
solutions of multivalued BDSDE \eqref{mbdsde:1}. Denote
\begin{eqnarray}
(\Delta Y(t), \Delta Z(t), \Delta U(t)):=(Y^1(t)-Y^2(t),
Z^1(t)-Z^2(t), U^1(t)-U^2(t)). \nonumber
\end{eqnarray}
By It\^{o}'s formula, we have
\begin{eqnarray}
&& e^{\beta t}|\Delta Y(t)|^2+\beta\int_t^Te^{\beta s}|\Delta
Y(s)|^2ds+\int_t^Te^{\beta s} |\Delta Z(t) |^2ds
+2\int_t^Te^{\beta s}\langle\Delta Y(s),\Delta U(s)\rangle ds
 \nonumber\\&=&
2\int_t^Te^{\beta s}\langle\Delta Y(s), (f(s,Y^1(s),Z^1(s),Y^1_s,
Z^1_s)-f(s,Y^2(s),Z^2(s),Y^2_s, Z^2_s))\rangle ds \nonumber\\&&
+\int_t^Te^{\beta s} |g(s,Y^1(s),Z^1(s),Y^1_s,
Z^1_s)-g(s,Y^2(s),Z^2(s),Y^2_s, Z^2_s)|^2ds \nonumber\\&&
+\int_t^Te^{\beta s}\langle\Delta Y(s),(g(s,Y^1(s),Z^1(s),Y^1_s,
Z^1_s)-g(s,Y^2(s),Z^2(s),Y^2_s, Z^2_s))dB(s)\rangle \nonumber\\&&
-2\int_t^Te^{\beta s}\langle \Delta Y(s),\Delta Z(s)dW(s)\rangle.
\nonumber
\end{eqnarray}
Since
\begin{eqnarray}
 \langle \Delta Y(s), \Delta U(s)\rangle\geq 0,\ \
d\mathbb{P}\times dt-a.e.\nonumber
\end{eqnarray}
Thus, as the same procedure as Lemma \ref{lemma4:3}, we can derive
the uniqueness of the solution. The proof is complete.
\\\\
\noindent{\bf Acknowledgments.} The first author is  grateful to
Prof. L. Maticiuc and Prof. A. R\c{a}\u{s}canu  for providing the
paper \cite{Diomande2013}.

\end{document}